\documentclass{amsart}

\newcommand{\R}{\mathbb{{R}}}
\newcommand{\N}{\mathbb{{N}}}

\newtheorem{theorem}{Theorem}[section]
\newtheorem{proposition}[theorem]{Proposition}

\theoremstyle{definition}

\theoremstyle{remark}
\newtheorem{remark}[theorem]{Remark}

\numberwithin{equation}{section}
\newcommand{\abs}[1]{\lvert#1\rvert}

\begin{document}

\title[Inequalities for the Gamma funtion and sections of $B_p^n$]
 {Inequalities for the Gamma function and estimates for
 the volume of sections of $B^n_p$}


\author[J. Bastero]{Jes\'us Bastero}
\address{Departamento de Matem\'{a}ticas, Facultad de Ciencias, Universidad de Zaragoza,
50009 Zaragoza, Spain}\email[(Je\'us
Bastero)]{bastero@posta.unizar.es}
\thanks{The first, the third and the forth authors were supported in part
by a DGES Grant (Spain).}

\author[F. Galve]{Fernando Galve}

\author[A. Pe\~{n}a]{Ana Pe\~{n}a}
\email[(Ana Pe\~{n}a)]{anap@posta.unizar.es}

\author[M. Romance]{Miguel Romance}
\thanks{The forth author was also supported by a FPI Grant (Spain).}
\email[(M. Romance)]{mromance@posta.unizar.es}


\date{March 28, 2000.}


\begin{abstract}

Let $B^n_p=\{(x_i)\in\R^n;\sum_1^n|x_i|^p\leq1\}$ and let $E$ be a
$k$-dimensional subspace of $\R^n$. We prove that $|E\cap
B^n_p|_k^{1/k}\geq  |B^n_p|_n^{1/n}$, for $1\leq k\leq (n-1)/2$
and $k=n-1$ whenever $1<p<2$. We also consider  $0<p<1$ and other
related cases. We obtain sharp inequalities involving Gamma
function in order to get these results.
\end{abstract}

\maketitle

\section*{Introduction}
In \cite{V} J.D. Vaaler proved that all $k$-dimensional sections
of the unit cube $[-1/2,1/2]^n$ in $\R^n$ have $k$-dimensional
volume bigger than or equal to $1$. If we write $B^n_p=\{x\in
\R^n; |x_1|^p+\cdots+ |x_n|^p\leq1\}$, $0<p< \infty$, M. Meyer
and A. Pajor extended Vaaler's result to the range $p\in\{1\}\cup
[2,\infty]$  in \cite {M-P}. More precisely, they showed that
\begin{equation}\label{E:Seccion}
|E\cap B^n_p|_k^{1/k}\geq  |B^n_p|_n^{1/n},
\end{equation}
 for all $E$, $k$-dimensional subspace of $\R^n$. K. Ball (see  \cite{B1},
\cite{B2}), using Brascamp and Lieb inequality, established
(\ref{E:Seccion}) for the $1$-dimensional sections of any ball in
$\R^n$ having a multiple of the euclidean ball as the ellipsoid of
maximal volume contained in it. More recently Schmuckenschl\"ager
(see \cite {S}), estimated the volume of the $(n-1)$-dimensional
sections of $B^n_p$, for $1<p<2$, but the proof of the inequality
he proposed was not correct. The aim of this paper is to give a
proof of the inequality appearing in \cite {S} (see proposition
\ref{P:Galve} below), to prove (\ref{E:Seccion}) for $1\leq k\leq
(n-1)/2$ and $1<p<2$ (see proposition \ref{T:1-p-2}) and to give
the right estimate for (\ref{E:Seccion}) for $0<p<1$ (see
proposition \ref{T:0-p-1}). Moreover, we prove (\ref{E:Seccion})
for $B=B^n_2\bigoplus_pB^n_2$, $1\leq p\leq2$, (see proposition
\ref{P:Suma}), answering a question raised to the authors by M.
Meyer.

In order to do this we need to establish sharp inequalities
involving the Gamma function which have their own interest. We
state and prove these inequalities in section 1 and the
corresponding estimates for the volume of sections are given in
section 2.

As usual we denote by $\Vert x\Vert_p=\left(\sum_1^n\abs
{x_i}^p\right)^{1/p}$, for $x=(x_1\ldots,x_n)\in \R^n$ and
$0<p<\infty$. If $A\subseteq \R^k$, $\abs{A}_k$ will represent the
$k$-dimensional Lebesgue measure in $\R^k$.

\section{Some inequalities for the Gamma function}

Throughout this paper we are going to use Stirling's formula, i.e.
$$\Gamma(1+x)=x^xe^{-x}\sqrt{2\pi x}\exp{\mu(x)}\qquad (x>0),$$
where $\mu (x)$ is a non increasing and non negative function for
$x\ge 1$ defined by $$\mu
(x)=\frac{1}{12x}-{1\over3}\int_0^{\infty}{p_3(t)\over(t+x)^3}
dt$$ where $p_3(t)$ is a $1$-periodic function that for
$t\in[0,1]$ is defined by $p_3(t)=t^3-{\frac{3}{2}}t^2+
{\frac{1}{2}}t$ (see, for instance \cite {Rem}, pg. 62). Notice
that $|p_3(t)|\le {1\over 20}$ for all $t\ge 0$.

\begin{proposition}\label{P:Acotacion} The following inequalities
hold:
\begin{enumerate}
\item $\displaystyle
\Gamma(1+x)^{2/x} \le {1\over 6}(x+1)(x+2),\qquad (x\ge 2)$
\item $\displaystyle
 \Gamma(1+x)^{2/x}\le {\frac{4}{23}}(x+1)(x+2),
\qquad (1\le x\le 2)$
\item $\displaystyle\Gamma(1+x)^{2/x} \ge {2\over
{\pi\sqrt{23}}}x\sqrt{(x+3)(x+7),}\qquad (x\ge 5)$
\item $\displaystyle\Gamma(1+x)^{2/x} \ge {1\over
e^2}(x+1)(x+2),\qquad (x\ge 1)$.
\end{enumerate}
\end{proposition}
\begin{proof}
We are only going to prove (1) and (3) because (2) and (4) can be
shown in a similar way.

\noindent (1) It is enough to prove that for every $x\ge 2$
\[f(x)=\log \Gamma(1+x) -{\frac{x}{2}}\log\left({\frac{1}
{6}}(x+1)(x+2)\right)\le 0.
\]
Let us compute $f^\prime$. $$ f'(x)=-{\frac{2-\log 6} {2}}+
\psi(x+1)-{\frac{\log {(x+1)(x+2)}}{2}}
 +{\frac{1}{x+2}}+{\frac{1/2}{x+1}}
$$ where $\psi(1+x)=\left(\log\Gamma(1+x)\right)'$. Now using that
\[\psi(1+x)<\log(x+1)-{\frac{1/2}{x+1}}-{\frac{1/12}{(x+1)^2}}
+{\frac{1/120}{(x+1)^4}}\] for $1+x>0$,  (see for instance
\cite{F}, section 541) and considering $y=\frac{1}{x+1}$ we
obtain that
\[\sup_{x\in [2,+\infty)} {f'(x)}\le \sup_{y\in (0,1/3]}
{g(y)},\] where
\[ g(y)=-{\frac{2-\log 6}{2}}-{\frac{1}{2}}\log
(y+1)+{\frac{y}{y+1}}-{\frac{y^2}{12}}+{\frac{y^4}{120}}.\] Since
$g$ is concave on $(0,1/3]$ and $g'(1/3)>0$, we get
\[ \sup_{y\in (0,1/3]}
{g(y)=g(1/3)<0}.\] Hence $f$ is a non increasing function on
$[2,+\infty)$ and so for every $x\ge 2$
\[\log \Gamma(1+x) -{\frac{x}{2}}\log\left({\frac{1}
{6}}(x+1)(x+2)\right)\le f(2)=0.\]

\noindent (3) Consider the function
$F:[5,+\infty)\longrightarrow\R$ defined by $$F(x)={4\over
x}\log\Gamma(x+1)-2\log x-\log (x+3) -\log (x+7)-
\log{4\over{23\pi^2}}.$$ We are going to show that $F(x)\ge 0$.
We have
\begin{align}
\notag x&F(x)= {4}\log\Gamma(x+1)-2x\log x-x\log (x+3) -x\log
(x+7)- x\log{4\over{23\pi^2}}\\ \notag\ge & {4}\log\left(x^x
e^{-x}\sqrt{2\pi x}\right)-2x\log x-x\log (x+3) -x\log (x+7)-
x\log{4\over{23\pi^2}}=G(x)\\ \notag
\end{align}
If we denote $\beta=-\log{{4e^4}\over{{23}\pi^2}}>0$, we get that
\begin{align}
\notag G'(x)= \left[\log\left(1-{3\over{x+3}}\right)+{3\over {5x}}
+{3\over{x+3}}+{\beta\over 2}\right]+
\left[\log\left(1-{7\over{x+7}}\right)+{7\over {5x}}
+{7\over{x+7}}+{\beta\over 2}\right].
\\
\notag
\end {align}
We have $$ h(y)= \log(1-y)+{{y}\over{5(1-y)}}+y+{\beta\over 2}>0$$
whenever $y\in[0,7/12]$, since $h$ is concave and
$h(0),h(7/12)>0$. It is clear that for every $x\ge 5$, $7/(x+7)$
and $3/(x+3)$ belong to the interval $[0,7/12]$, therefore
$G^\prime (x)>0$ on $[5,+\infty)$. Hence we conclude that $G(x)\ge
G(5)>0$ for all $x\in [5,+\infty)$, and so $F(x)>0$, for all
$x\ge 5$.
\end{proof}

\begin{proposition}\label{P:Galve} Let $1/2\le x\le 1$ and $y\ge 2$, then
\begin{equation}\label{E:Gamma}
{ \Gamma (1 + {xy})^{1 + {2 \over y}} \over \Gamma \left(1 +
{(y+2)x}\right)} \leq { \Gamma (1 + {x})^3 \over \Gamma (1 + {3x})
}
\end{equation}
\end{proposition}

\begin{proof} First of all note that for every $y\ge 2$, (\ref{E:Gamma})
holds for $x=1$ and $x=1/2$, simply using proposition
\ref{P:Acotacion} and because $${{\Gamma(1+y)^{1+2/y}}\over
{\Gamma(3+y)}}= {{\Gamma(1+y)^{2/y}}\over (y+2)(y+1)}\le {1\over
6}$$ and $${{\Gamma(1+{y\over 2})^{1+2/y}}\over {\Gamma(2+{y\over
2})}}= {{\Gamma(1+{y\over 2})^{2/y}}\over {1+{y\over 2}}}\le
\frac2{\sqrt{23}}\sqrt{\frac{y+4}{y+2}}\leq \frac{\pi}6.$$ Now we
only have to prove that for every $y\ge 2$ the function $f_y :
[{1 \over 2},1] \to {\bf R}$ defined by $f_y(x) = 3
\log\Gamma(1+x) - \log\Gamma(1+3x) - (1+{2 \over
y})\log\Gamma(1+xy) + \log\Gamma(1+(y+2)x) $ is concave. If we
compute its derivate, we obtain $$ f_y''(x) = 3 \psi'(1+x) - 9
\psi'(1+3 x) - y(y+2) \psi'(1+xy)+ (y+2)^2 \psi'(1+(y+2) x) $$
Next we use that there exists a function
$\theta:(0,+\infty)\longrightarrow [0,1]$ such that for every
$z>0$ $$\psi'(z)={1\over z}+{1\over {2z^2}}+{1\over
{6z^3}}-{{\theta(z)}\over {30z^5}}$$ (see \cite{F} or \cite{G})
then we get that
\begin{multline}
\notag f_y''(x) = \underbrace{-{1 \over {(x+1)(x+{1 \over 3})}}
\left[ 1 + {4 \over 3} \cdot { {x + {1 \over 2}} \over {(x+1)(x+{1
\over 3})}} \right]}_{\bf S1} + \underbrace{{-1 + 9x^2 + 12x^3
\over (1+x)^3 \, (1+3x)^3}}_{\bf S2} \\ + \underbrace{{ 1/y \over
{(x+{1 \over y})(x+{1 \over y+2})} } \left[ 1 + {2(y+1) \over
y(y+2)} \cdot { {x + {1 \over y+1}} \over {(x+{1 \over y})(x+{1
\over y+2})}} \right]}_{\bf S3} + \alpha(x,y), \\
\end{multline}
where \begin{align} \notag\alpha(x,y)&\le {9\over
{30(1+3x)^5}}+{{y(y +2)}\over{30(1+xy)^5}}- {y (y + 2) \over 6 (1
+ xy)^3} + {(y + 2)^2 \over 6 (1 + (y + 2) x)^3}\\ &\le {9\over
{30(1+3x)^5}}+{{y(y +2)}\over{30(1+xy)^5}}\notag
\end{align}

\noindent Notice that {\bf S1} and {\bf S3} can be deduced from
the identity
\begin{multline}
\notag {-k(k+2) \over (1 + k x)} + {-k(k+2) \over 2(1 + k x)^2} +
{(k+2)^2 \over 1 + (k + 2) x} + {(k+2)^2 \over 2(1 + (k + 2) x)^2}
\\
={{1 \over k} \over (x + {1 \over k})(x + {1 \over k+2})} \left[ 1
+ {2(k+1) \over k(k+2)} \cdot { {x + {1 \over k+1}} \over {(x+{1
\over k})(x+{1 \over k+2})}} \right]  \\
\end{multline}
applied to $k=1$ and $k=y$ respectively.

Now we are going to study each summand separately:

\noindent {\bf S1}: It is easy to check that $$ \min_{x \in [{1
\over 2},1]} \left( 1 + {4 \over 3} \cdot { x + {1 \over 2} \over
(x + 1)(x + {1 \over 3}) } \right) = {7 \over 4} \qquad (x=1)$$

\noindent {\bf S2}: It can be shown that $$ \max_{x \in [{1 \over
2},1]} {9 x^2 + 12 x^3 -1 \over (1 + x)^2 \, (1 + 3 x)^2}< \max_{x
\in [{1 \over 2},1]} {9 x^2 + 12 x^3 \over (1 + x)^2 \, (1 + 3
x)^2}
=
{21 \over 64} \qquad (x=1)$$

\noindent {\bf S3}: $ \displaystyle \max_{x \in [{1 \over 2},1]}
\left( 1 + {2(y+1) \over y(y+2)} \cdot { {x + {1 \over y+1}} \over
{(x+{1 \over y})(x+{1 \over y+2})}} \right) = 1 + {4 (y + 3)
\over (y + 2)(y + 4)}\qquad (x=1/2)$. Notice that since $$ \max_{x
\in [{1 \over 2},1]} { (x + 1)(x + {1 \over 3}) \over (x + {1
\over y})(x + {1 \over y +2}) } = {5 y \over y + 4} \qquad
(x=1/2)$$ we get that for every $1/2\le x\le 1$ and every $y\ge 2$
\begin{align} \notag { 1/y \over {(x+{1 \over y})(x+{1 \over y+2})}
} \left[ 1 + {2(y+1) \over y(y+2)} \cdot { {x + {1 \over y+1}}
\over {(x+{1 \over y})(x+{1 \over y+2})}} \right]  &\le { {5 \over
y + 4} \cdot \left(1 + {4 (y + 3) \over (y + 2)(y + 4)}\right)
\over (x + 1)(x + {1 \over 3}) } \\ \notag &\le {{55/
36}\over{(x+1)(x+1/3)}}.\\ \notag
\end{align}
On the other hand $$\alpha(x,y)\le   {9 \over {30(1+{3/2})^5}}
 +  { {y(y+2)} \over {30 \left( y+2 \over
2\right)^5} }\le {4277\over 375000}(\sim 0.012).$$ Therefore  we
obtain that for each $y\ge 2$,
\begin{align}\notag
f_y''(x)&\le{\frac{-21/4}{(1+x)(1+3x)}}+{\frac{21/64}
{(1+x)(1+3x)}}+{\frac{55/36}{(x+1)(x+1/3)}}+{\frac{4277}{375000}}\\
& ={\frac{4277}{375000}}+ { {-{\frac{105}{64}} + {\frac{55}{36}}}
\over {(x + 1)(x + {1 \over 3}) }} < 0 \notag
\end{align}
 for all $1/2\le x\le 1$.
\end{proof}

\begin{proposition}\label{P:2Variables}
 The function
$g:[1/2,+\infty)\times[2,+\infty)\longrightarrow (0,+\infty)$
defined by
\[
g(x,y)={\frac{y}{{\Gamma(1+y)}^{1/y}}}{\frac{{\Gamma(1+xy)}^{1/y}}
{{y^x} \Gamma(1+x)}}
\]
verifies:

\begin{enumerate}
\item For every $y\ge 2$, $g(\cdot,y)$ is non
increasing in $[1/2,+\infty)$.

\item For every $x\ge 1$, $g(x,\cdot)$ is non
increasing in $[2,+\infty)$ and for every $1/2\le x\le 1$,
$g(x,\cdot)$ is non decreasing in $[2,+\infty)$.
\end{enumerate}
\end{proposition}

\begin{proof}
(1)  Let $y\ge 2$. By using Stirling's formula it is easy to see
that
\[ h(x,y)=g(x,y){\frac{\Gamma(1+y)^{1/y}}{y}}=y^{\frac{1}{2y}}(2\pi x)^
{\frac{1}{2y}-\frac {1} {2}}\exp{\left({\frac{\mu(xy)}
{y}}-\mu(x)\right)}. \]
 Since
$|p_3(t)|\leq1/20$, we have
\begin{align}
\notag \frac{\partial (\log h)}{\partial x} (x,y)=& -{{y-1}\over
{2xy}}+{\frac{1}{12x^2}}-{1\over{12x^2y^2}}+
{}{\int_0^{+\infty}\left({{p_3(t)}\over{(t+xy)^4}}-
{{p_3(t)}\over{(t+x)^4}}\right)dt}\\ \notag\le &-{{y-1}\over
{2xy}}+{1\over{12x^2}}-{1\over{12x^2y^2}}+
{}{1\over{20}}{\int_0^{+\infty}\left({1\over{(t+x)^4}}-
{1\over{(t+xy)^4}}\right)dt}\\ \notag \le &
{1\over{2x}}\left(1-{1\over y}\right)\left(-1+{1\over {4x}}+
{7\over {120x^2}}\right)<0
\end{align}
for all $x\ge 1/2$, (note that this result can be extended to $x$
strictly smaller than $1/2$). Therefore $\log h(\cdot,y)$ is a
non increasing function in $[1/2,+\infty)$ and so it is $g(x,y)$.

\noindent (2) Let $x\ge 1/2$. If we use again Stirling's
expression of the Gamma function we have $$g(x,y)=
{{e^{1-x}x^x}\over{\Gamma(1+x)}}x^{1/2y} \exp\left({1\over
y}(\mu(xy)-\mu(y))\right)$$ Consider the function $$\phi(x,y)=
{1\over 2y}\log x +{1\over y}(\mu(xy)-\mu(y))$$ defined for $y\ge
1$ and $x\ge 1/2$. $$\phi(x,y)= {1\over 2y}\log x +{1\over
12y^2}({1\over x}-1)-{1\over 3y}\int_0^\infty p_3(t)\left({1\over
(xy+t)^3}-{1\over (y+t)^3}\right)dt$$ Then $$ {\partial^2
\phi\over
\partial x\partial y}(x,y)= -{1\over {2x y^2}}+{1\over{6x^2 y^3}}-
4x\int_0^\infty {p_3(t)\over (xy+t)^5}dt. $$ Since $\max
\{|p_3(t)|; t\ge0\}\leq {1\over 20}$, we achieve
\begin{align}
\notag {\partial^2 \phi\over \partial x\partial y}(x,y) &\leq
-{1\over {2 y^2 x}}+{1\over{6x^2 y^3}}+{1\over{ 20x^3 y^4}}\notag
\\ \notag &\leq -{1\over {2y^2 x}}\left(1-{1\over {6x}}-{1\over
{40x^2}}\right)<0
\end{align}
Hence if $x\ge 1$, $$ {\partial \phi\over \partial
y}(x,y)\le{\partial \phi\over
\partial y}(1,y)=0$$ for all
$y\in [2,\infty)$ and if $1/2\le x\le 1$ $$0={\partial \phi\over
\partial y}(1,y)\le {\partial \phi\over \partial
y}(x,y)$$ for all $y\in [2,\infty)$ and thus, the result holds.
\end{proof}

\begin{proposition}\label{P:4gamas}
 The following inequality holds
\[\Gamma(1+2x)\Gamma\left(1+\frac x2\right)^2\geq
2^x\Gamma(1+x)^2\Gamma\left(1+\frac{2x-1}2\right)^{2x/(2x-1)}\]
for all $x\geq 5/2$.
\end{proposition}

\begin{proof} We apply Stirling formula and so, we only need to achieve
\begin{equation}\label{P:uno}
\left(\frac x{x-1/2}\right)^{x+1/2}\geq \sqrt 2
\left((2x-1)\pi\right)^{1/(4x-2)}
\end{equation}
and
\begin{equation}\label{P:dos}
2\mu(x)+\frac{2x}{2x-1}\mu\left(\frac{2x-1}2\right)-\mu(2x)
-2\mu\left(\frac x2\right)\leq0.
\end{equation}

The inequality (\ref{P:uno}) is deduced from the fact that the
function
\[F(y)=y(1+y/2)\log(1+y^{-1})-\frac y2\log2-\frac 12\log(\pi y)\] is
convex for $y>0$. In particular since $F'(4)>0$ and $F(4)>0$ we
deduce $F(y)>0$ for all $y\geq4$ and so the inequality is true
for $x\geq 5/2$ (consider $2x-1=y$).

In order to show (\ref{P:dos}) we use the corresponding expansion
and we have
\begin{align}
2&\mu(x)+\frac{2x}{2x-1}\mu\left(\frac{2x-1}2\right)-\mu(2x)
-2\mu\left(\frac x2\right) = -\frac5{24x}+\frac
x{3(2x-1)^2}\notag\\
&-\frac13\int_0^\infty p_3(t)\left(\frac2{(x+t)^3}
+\frac{2x}{2x-1}\frac1{(t+(2x-1)/2)^3}-\frac1{(2x+t)^3}
-\frac2{(t+x/2)^3}\right)dt\notag\\
& =-\frac5{24x}+\frac x{3(2x-1)^2}+\frac13\int_0^\infty
p_3(t)\left(\frac2{(t+x/2)^3}-\frac2{(x+t)^3}\right)\,dt \notag\\
& -\frac13\int_0^\infty
p_3(t)\left(\frac{2x}{2x-1}\frac1{(t+(2x-1)/2)^3}-\frac1{(2x+t)^3}\right)
\,dt \notag\\
 & \leq\frac1{24x}\left(-5+\frac{23}{20x}+\frac8{(2-1/x)^2}
+\frac{8/5}{(2-1/x)^2(2x-1)}\right) \leq -\frac{0.632}{24x}<0
\notag
\end{align}
since the function
\[-5+\frac{23}{20x}+\frac8{(2-1/x)^2}
+\frac{8/5}{(2-1/x)^2(2x-1)}\] is non increasing for $x\geq2$.
Therefore the result follows.
\end{proof}

\section{The volume of central sections of the
unit ball in ${\ell^n_p}$, $0< p< 2$}

We apply the preceding inequalities to estimate the volume of the
$k$-dimensional sections of $B^n_p$, stated in the introduction.

\begin{proposition}(see \cite{S}).
Let $n\in\N$, $n\geq 2$, $p\in [1,2]$ and let $E$ be any
$(n-1)$-dimensional subspace in $\R^n$. Then $$\left|E\cap
B_p^n\right|_{n-1}^{1/{n-1}}\ge \left|B_p^n\right|_n^{1/n}.$$
\end{proposition}

\begin{proof}
Let $E$ be a hyperplane in $\R^n$. A well known result (see
\cite{H}) ensures that

$$\left|B_p^n \cap E\right|_{n-1} L_{B_p^n}\ge{1\over {\sqrt
{12}}}\left| B_p^n\right|_n^{(n-1)/n}$$ where $L_{B_p^n}$ (the
isotropy constant) is $$ {L_{B_p^n}^2}={{\Gamma(1+ {3\over p} )
\Gamma(1+{n\over p})^{1+2/n}} \over{12 \Gamma(1+{{n+2}\over
p})\Gamma (1+{1\over p})^3}}$$ (see \cite{Mi-P}). Hence it is
enough to prove that $${{\Gamma(1+ {3\over p} ) \Gamma(1+{n\over
p})^{1+2/n}} \over{ \Gamma(1+{{n+2}\over p})\Gamma (1+{1\over
p})^3}} \le 1$$ for all $n\ge 2$ and all $1\le p\le 2$. Notice
that this follows from proposition \ref{P:Galve}.
\end{proof}

\begin{proposition}\label{T:1-p-2}
Let $n\in\N$, $n\geq 2$, $p\in [1,2]$ and let $E$ be any
$k$-dimensional subspace in $\R^n$ with $1\le k\le {{n-1}\over
2}$. Then $$\left|E\cap B_p^n\right|_k^{1/k}\ge
\left|B_p^n\right|_n^{1/n}.$$
\end{proposition}

\begin{proof}
Acording to the K. Ball's result quoted in the introduction,  we
only have to consider the case $n\ge 5$.

H\"older's inequality implies that $$\left|E\cap
B_p^n\right|_k^{1/k}\ge\left|E\cap n^{{1\over 2}-{1\over p}}
B_2^n\right|_k^{1/k}=n^{{1\over 2}-{1\over p}}
\left|B_2^k\right|_k^{1/k}$$ (in fact $n^{1/2-1/p}B^n_2$ is the
ellipsoid of maximal volume contained in $B^n_p$). Hence it is
enough to show that $$n^{{1\over 2}-{1\over
p}}\left|B_2^k\right|_k^{1/k}\ge \left|B_p^n\right|_n^{1/n}$$ for
all $1\le p\le 2$ and for all $1\le k\le {{n-1}\over 2}$, that is,
\begin{equation}\label{E:Segunda}
{{\Gamma({n\over p}+1)}^{1/n} \over{n^{1/p} \Gamma({1\over
p}+1)}}\ge {{\Gamma({k\over 2}+1)}^{1/k} \over{n^{1/2}
\Gamma({1\over 2}+1)}}
\end{equation}
(see for instance \cite{Pisier}).

By using proposition \ref{P:2Variables}, for every $1\le p\le 2$
we get that $g(1/p,n)\ge g(1,n)$, therefore it is enough to prove
(\ref{E:Segunda}) for $p=1$. Furthermore, since $\Gamma (1+x)$ is
log-convex on $[0,+\infty)$, the function
$f(x)=\Gamma(1+x)^{1\over x}$ is a non decreasing funtion on
$[0,+\infty)$, so
\[{\Gamma({k\over 2}+1)}^{1/k}
\le{\Gamma({{n-1}\over 4}+1)}^{2\over{n-1}}\] for all $1\le k\le
{{n-1}\over 2}$. On the other hand, since ${{n-1}\over 4}\ge 1$,
we can use proposition \ref{P:Acotacion}, $(1)$ and $(2)$, and we
obtain that $${{\Gamma({{n-1}\over 4}+1)}^{2\over{n-1}} \over{
\Gamma({1\over 2}+1)}}\le {2\over {\sqrt{\pi}}} \left({1\over
92}(n+3)(n+7)\right)^{1/4}.$$ Thus, it sufficies to show that
\begin{equation}
\Gamma(1+n)^{4/n}\ge{4\over{23\pi^2}}n^2(n+3)(n+7)
\end{equation}
for all integer $n\ge 5$ and this is proposition
\ref{P:Acotacion}.
\end {proof}

\begin{remark} If we consider
$$K=\left\{(x_1,\ldots,x_m)\in \R^n \times\dots\times\R^n;
\enspace \|x_1\|_2^p+ \dots + \|x_m\|_2^p \le 1 \right\},$$ with
$1\le p\le 2$ and $n,m\in \N$, and  we use the same method as in
proposition \ref{T:1-p-2}, it can be shown that for every
$k$-dimensional linear subspace $E$ in $\R^{nm}$ with $1\le k\le
{{nm-1}\over 2}$
\begin{equation}\label{E:Suma2} \left|E\cap
K\right|_k^{1/k}\ge \left|K\right|_{nm}^{1/{nm}}
\end{equation}
for all $p\in [1,2]$ and all $n,m\in \N$. The only new tool we
need is the inequality
\begin{equation}
{{n^{1/2}\Gamma(1+{n\over 2})^{1/n}}\over{\Gamma (1+{1\over 2})
\Gamma(n+1)^{1/n} }}\ge 1 \qquad (n\ge 1)
\end{equation}
which is a consecuence of proposition \ref{P:2Variables}. Moreover
we can achieve the inequality (\ref{E:Suma2}) for all $1\leq k\leq
2n$, when $m=2$, and this way extends the results in \cite{M-P} in
this case, as it is shown in the following result.
\end{remark}

\begin{proposition}\label{P:Suma}
Let $1\le p\le 2$, $n\in\N$ and $$K=\left\{(x_1,x_2)\in \R^n
\times\R^n; \enspace \|x_1\|_2^p+ \|x_2\|_2^p \le 1 \right\}$$
then (\ref{E:Suma2}) holds for all $k$-dimensional subspace in
$\R^{2n}$, with $1\le k\le 2n$.
\end{proposition}
\begin{proof}
Following the same methods than proposition \ref{T:1-p-2}, we
only have to prove
\[\Gamma(1+2n)\Gamma(1+\frac n2)^2\geq
2^n\Gamma(1+n)^2\Gamma(1+\frac{2n-1}2)^{2n/(2n-1)}\] for $n\geq
2$. The case $n\geq 3$ is proposition \ref{P:4gamas} and $n=2$ can
be checked directly.
\end{proof}

Next we are going to estimate the volume of the sections through
the origin for the $p$-balls $B^n_p$, $0<p<1$. We should notice
that Koldobsky (see \cite{K}) studied  this problem for the
particular case of central hyperplane sections. He computed the
volume of these sections in terms of the Fourier transform of a
power of the radial function, for every $p>0$, and he applied
this result to confirm the conjecture of Meyer and Pajor on the
minimal volume of these particular sections of the unit $p$-balls
$B^n_p$, $0<p<2$.

\begin {proposition}\label{T:0-p-1} Let $E$ be any  $k$-dimensional
subspace of $\R^n$, $1\leq k\leq n$ and let $0<p<1$. Then
\[ |{B^n_p\cap E}|_k^{1/k}\ge {\frac{e^{1-1/p}}{\Gamma(1+1/p)p^{1/p}}}
 |{B^n_p}|_n^{1/n}\]
 and the constant
$${\frac{e^{1-1/p}}{\Gamma(1+1/p)p^{1/p}}}\in (0,1)$$ is the good
order of magnitude for fixed $n$ when $p\longrightarrow 0^+$.
\end {proposition}

\begin {proof}
We use   the results from \cite{M-P}
\begin{align}
\notag   |{E\cap B^n_p}|_k^{1/k}&\geq n^{1-1/p}
|{E\cap B^n_1}|_k^{1/k}\\
\notag &\geq  n^{1-1/p}  |{B^n_1}|_n^{1/n}\\
\notag &=n^{1-1/p}{\Gamma(1+n/p)^{1/n}\over (n!)^{1/n}\Gamma(1+1/p)} |{
B^n_p}|_n^{1/n}.\notag
\end{align}
By proposition \ref{P:2Variables}
$$n^{1-1/p}{\Gamma(1+n/p)^{1/n}\over (n!)^{1/n}\Gamma(1+1/p)}$$ is
non increasing with $n$ and this implies the result, since
$$\lim_{n\to \infty}n^{1-1/p}{\Gamma(1+n/p)^{1/n}\over (n!)^{1/n}
\Gamma(1+1/p)}={e^{1-1/p}\over\Gamma(1+1/p)p^{1/p}}.$$

\noindent Note that this value belongs to $(0,1)$. Indeed
\begin{align}
\notag {e^{1-1/p}\over\Gamma(1+1/p)p^{1/p}}&={e\sqrt p\over
\sqrt{2\pi}}\exp  (-\mu(1/p))\\ \notag &={e\over \sqrt{2\pi}}\exp
\left({1\over 2}\log p-{p\over 12}+ {1\over
3}\int_0^\infty{p_3(t)\over (1/p+t)^3}dt\right)
\end {align}
and
\begin{align}
\notag {d\over dp}\left({1\over 2}\log p-{p\over 12}+ {1\over
3}\int_0^\infty{p_3(t)\over (1/p+t)^3}dt\right)&= {1\over 2p}
-{1\over 12}+{1\over p^2}\int_0^\infty{p_3(t)\over (1/p+t)^4}dt\\
\notag &\geq {1\over 2p} -{1\over 12}-{p\over 60}>0.
\end{align}

\noindent Finally we show that the result is sharp. It is easy to
check that $${e^{1-1/p}\over\Gamma(1+1/p)p^{1/p}}\sim e{\sqrt
{\frac {p}{2\pi}}},\qquad {\rm when}\enspace p\longrightarrow
0^+$$ and if we consider the 1-dimensional subspace $E_0={\rm
span\,}\{(1,\ldots,1)\}\subseteq \R^n$ then it is easy to prove
that $${{\left|B_p^n\cap E_0\right|_1}\over
{\left|B_p^n\right|_n^{1/n}}}=n^{1/2-1/p}{\Gamma(1+n/p)^{1/n}\over
\Gamma(1+1/p)}\sim \frac{p^{{1\over 2}-{1\over {2n}}}n^{{1\over
2}+{1\over{2n}}} }{(2\pi)^{{1\over 2}-{1\over 2n}}}\qquad
(p\longrightarrow 0^+).$$
\end {proof}

\begin{remark} If we now consider
$$K=\left\{(x_1,\ldots,x_m)\in \R^n \times\dots\times\R^n;
\enspace \|x_1\|_1^p+ \dots + \|x_m\|_1^p \le 1 \right\},$$ with
$0< p\le 1$ and $n,m\in \N$, and we use the same ideas as in
proposition \ref{T:0-p-1} it can be shown that for every
$k$-dimensional linear subspace $E$ in $\R^{nm}$, $1\leq k\leq
nm$, $$\left|E\cap K\right|_k^{1/k}\ge
{\frac{e^{1-1/p}}{\Gamma(1+1/p)p^{1/p}}}
\left|K\right|_{nm}^{1/{nm}}$$ for all $p\in (0,1]$ and all
$n,m\in \N$.
\end{remark}

\bibliographystyle{amsplain}

\end{document}